\documentclass[11pt]{amsart}

\usepackage{algorithm}
\usepackage{algpseudocode,lineno}
\usepackage{nicefrac}

\usepackage{amssymb,amscd,amsthm,verbatim,amsmath,color,fancyhdr,mathrsfs}
\usepackage{graphicx}
\usepackage{turnstile}

\usepackage[top=.95in, bottom = 1 in, left=1in, right = .6in]{geometry}
\usepackage{styfile} 

\newtheorem{theorem}{Theorem}[section]

\newtheorem{lemma}[theorem]{Lemma}

\theoremstyle{definition}

\newtheorem{conjecture}{Conjecture}

\DeclarePairedDelimiterX{\bracket}[3]{#1}{#2}{#3}
\newcommand{\round}[1]{\bracket*{(}{)}{#1}}

\providecommand{\newoperator}[3]{\newcommand*{#1}{\mathop{#2}#3}}
\providecommand{\renewoperator}[3]{\renewcommand*{#1}{\mathop{#2}#3}}

\renewoperator{\Re}{\mathrm{Re}}{\nolimits}
\renewoperator{\Im}{\mathrm{Im}}{\nolimits}

\DeclarePairedDelimiterXPP{\nrm}[2]{}{\lVert}{\rVert}{\ensuremath{_{#1}}}{\ifblank{#2}{\:\cdot\:}{#2}}

\newcommand{\abs}[1]{\bracket*{\lvert}{\rvert}{#1}}

\DeclarePairedDelimiterXPP\prob[1]{\mathbb{P}}{\lbrace}{\rbrace}{}{#1} 

\DeclarePairedDelimiterXPP\probability[2]{\mathbb{P}_{#1}}{\lbrace}{\rbrace}{}{#2} 

\DeclarePairedDelimiterXPP\expectation[1]{\mathbb{E}}{\lbrack}{\rbrack}{}{#1} 

\DeclarePairedDelimiterXPP\expectationdist[2]{\mathbb{E}_{#1}}{\lbrack}{\rbrack}{}{#2} 

\DeclarePairedDelimiterXPP\variance[1]{\mathrm{Var}}{\lbrack}{\rbrack}{}{#1} 

\DeclarePairedDelimiterXPP\variancedist[2]{\mathrm{Var}_{#1}}{\lbrack}{\rbrack}{}{#2} 

\DeclarePairedDelimiterXPP\covariance[2]{\mathrm{Cov}}{(}{)}{}{#1,\mathopen{}#2} 

\DeclarePairedDelimiterXPP\covariancedist[3]{\mathrm{Cov}_{#1}}{(}{)}{}{#2,\mathopen{}#3} 

\newoperator{\supp}{\mathrm{supp}}{\nolimits}

\makeatletter
\providecommand*{\diff}
{\@ifnextchar^{\DIfF}{\DIfF^{}}}
\def\DIfF^#1{
	\mathop{\mathrm{\mathstrut d}}
	\nolimits^{#1}\gobblespace}
\def\gobblespace{
	\futurelet\diffarg\opspace}
\def\opspace{
	\let\DiffSpace\!
	\ifx\diffarg(
	\let\DiffSpace\relax
	\else
	\ifx\diffarg[
	\let\DiffSpace\relax
	\else
	\ifx\diffarg\{
	\let\DiffSpace\relax
	\fi\fi\fi\DiffSpace}

\providecommand*{\pdiff}
{\@ifnextchar^{\pDIfF}{\pDIfF^{}}}
\def\pDIfF^#1{
	\mathop{\mathrm{\mathstrut \partial}}
	\nolimits^{#1}\gobblespace}
\def\gobblespace{
	\futurelet\diffarg\opspace}
\def\opspace{
	\let\DiffSpace\!
	\ifx\diffarg(
	\let\DiffSpace\relax
	\else
	\ifx\diffarg[
	\let\DiffSpace\relax
	\else
	\ifx\diffarg\{
	\let\DiffSpace\relax
	\fi\fi\fi\DiffSpace}

\makeatother

\DeclarePairedDelimiterX\Set[1]\{\}{
	
	#1
}

\newcommand{\N}{\mathbb{N}}

\newcommand{\Ccal}{\mathcal{C}}
\newcommand{\Dcal}{\mathcal{D}}

\newcommand{\Ocal}{\mathcal{O}}

\newcommand{\eq}{\begin{equation}}
\newcommand{\en}{\end{equation}}

\newcommand{\Sym}[1]{\mathfrak{S}_{#1}}

\newoperator{\tensor}{\otimes}{}

\title{On log-concavity of the number of orbits in commuting tuples of permutations}

\author{Raghavendra Tripathi}
\address{Raghavendra Tripathi\\ Department of Mathematics \\ University of Washington\\ Seattle WA 98195, USA\\ {Email: raghavendra.cic@gmail.com}}

	

\date{\today}

\begin{document}

\begin{abstract}
Denote by $A(p, n, k)$ the number of commuting $p$-tuples of permutations on $[n]$ that have exactly $k$ distinct orbits. It was conjectured in~\cite{abdesselam2023log} that $A(p, n, k)$ is log-concave with respect to $k$ for every $p\geq 2, n\geq 3$, and the log-concavity was proved in ``$p=\infty$" case. In this paper, we prove that for $k=n-\alpha$, the log-concavity for $A(p, n, k)$ holds for every $p\geq 2$ for sufficiently large $n$.
\end{abstract}

\maketitle



\section{Introduction}\label{sec:Intro}
For $n\in\N$, let $\Sym{n}$ denote the symmetric group on the set $[n]=\{1, \ldots, n\}$. For any permutation $\sigma\in \Sym{n}$, let $\Ocal(\sigma)$ denote the set of orbits under the action of $\sigma$ on $[n]$. For any $1\leq k\leq n$, let us define 
\[ A(1, n, k) = \{\sigma:\Sym{n}: \abs{\Ocal(\sigma)}=k\}\;\;.\]
That is, $A(1, n, k)$ is the number of permutations $\sigma\in \Sym{n}$ that has exactly $k$ orbits. It is well-known that $A(1, n, k)$ equals the unsigned Stirling number of the first kind, and using this it can be shown (see, for example,~\cite{abdesselam2023local}) that the $A(1, n, k)$ is log-concave, that is, it satisfies 
\[ A(1, n, k)^2\geq A(1, n, k-1)A(1, n, k+1), \qquad 2\leq k\leq n-1, \quad n\geq 3.\]

The quantity $A(1, n, k)$ was generalized to $A(p, n, k)$ as the number of pairwise commuting $p$-tuples of permutations such that the subgroup generated by them has exactly $k$ orbits. More precisely, let $p, n\in \N$ and let $\Ccal_{p, n}$ be the set of commuting $p$-tuples of permutations on $\Sym{n}$. That is, 
\[\Ccal_{p, n}= \{(\sigma_1, \ldots, \sigma_p)\in (\Sym{n})^p: \sigma_i\sigma_j=\sigma_j\sigma_i \;\;\;\forall i, j\in [p]\}\;.\]
Furthermore, for $\mathbf{\sigma}=(\sigma_1, \ldots, \sigma_p)\in \Ccal_{p, n}$, let $\langle \mathbf{\sigma}\rangle$ be the subgroup of $\Sym{n}$ generated by $\sigma_1, \ldots, \sigma_p$. Let $\Ocal(\mathbf{\sigma})$ denote the number of orbits of $[n]$ under the action of $\langle \mathbf{\sigma}\rangle$. For $1\leq k\leq n$, define $\Ccal_{p, n, k}$ to be the set 
\[\Ccal_{p, n, k} = \{\mathbf{\sigma}\in C_{p, n}: \abs{\Ocal(\mathbf{\sigma})}=k\}\;,\]
and let $A(p, n, k) = \abs{\Ccal_{p, n, k}}$. It was conjectured in~\cite{abdesselam2023log} that $A(p, n, k)$ is log-concave for all $p$ and $n$. More precisely, 
\begin{conjecture}[~\cite{abdesselam2023log}]
\label{conj:Conjecture}
Let $p\geq 2$ and let $n\geq 3$. Then, 
\[ A(p, n, k)^2\geq A(p, n, k-1)A(p, n, k+1)\;,\qquad 2\leq k\leq n-1\;.\]
\end{conjecture}

The case $p=2$ was already conjectured in~\cite{heim2021horizontal} and it implies the so-called unimodality conjecture for Nekrasov-Okunkov polynomials. We refer the reader to~\cite{ abdesselam2023log, abdesselam2024bijection} for more in-depth discussion about this conjecture and~\cite{heim2019conjectures, heim2020formulas,hong2021towards,zhang2022log} for more discussion as well as the work focused primarily on the case $p=2$. 

The quantity $A(p, n, k)$ admits a neat and clear formula given by
\begin{equation}
    \label{eqn:masterEq}
    A(p, n, k) = \frac{n!}{k!} \times \sum_{n_1, \ldots, n_k\geq 1} \mathbbm{1}\{n_1+\cdots+n_k=n\}\times \prod_{i=1}^{k} \frac{B(p, n_i)}{n_i}\,
\end{equation}
where $\mathbbm{1}\{n_1+\cdots+n_k=n\}=1$ if the tuple $(n_1, \ldots, n_k)$ satisfies $n_1+n_2+\cdots+n_k=n$ and $0$ otherwise, and $B(p, m)$ is a multiplicative function defined as
\begin{equation}
    B(p, m) = \sum_{s_1|s_2|\cdots|s_{p-1}|m} s_1s_2\ldots s_{p-1}\;.
\end{equation}
We used the notation $x|y$ to denote $x$ divides $y$. While the formula~\eqref{eqn:masterEq} can be obtained from the results in~\cite{bryan1998orbifold}, the authors in~\cite{abdesselam2024bijection} give a bijective proof of~\eqref{eqn:masterEq} in the spirit of enumerative combinatorics. We highly recommend an interested reader to go through the introduction in~\cite{abdesselam2024bijection,westbury2006universal,nekrasov2006seiberg}, and references therein for more details. 

Conjecture~\ref{conj:Conjecture} for $p=2$ and $k=n-1$ was proved in~\cite[Proposition 4]{heim2020formulas}. Using the explicit values for $A(p, n, k)$ for $k=n-1, n-2, n$, thanks to~\eqref{eqn:masterEq}, the authors in~\cite[Proposition 3.1]{abdesselam2024bijection} also resolve the Conjecture~\ref{conj:Conjecture} for the case $k=n-1$ for all $p\geq 2, n\geq 3$ and also verify the Conjecture~\ref{conj:Conjecture} using Mathematica for $p=3, 4, 5$ and $n\leq 100$. In~\cite[Theorem 1.1]{abdesselam2023log}, the Conjecture~\ref{conj:Conjecture} is proved for the case $p=\infty$, that is, it is shown that for all $n\geq 3$ and all $2\leq k\leq n-1$ the following inequality holds:
\[ \liminf\limits_{p\rightarrow \infty}\frac{A(p, n, k)^2}{A(p, n, k-1)A(p, n, k+1)}\geq 1\;.\]

In this paper, we complement the above result in another direction. We show that the Conjecture~\ref{conj:Conjecture} holds for all $k$ close to $n$ when $n$ is sufficiently large. More precisely, we prove the following. 

\begin{theorem}
\label{thm:MainTheorem}
Let $p\geq 2$ and let $\alpha\in \N$ be fixed. Let $k=n-\alpha$. Then, there exists $N=N(\alpha, p)$ depending only on $p$ and $\alpha$ such that 
\[ A(p,n, k)^2\geq A(p, n, k-1)A(p, n, k+1)\;,\]
for all $n\geq N$. 
\end{theorem}

The proof of Theorem~\ref{thm:MainTheorem} is based on analyzing the first-order asymptotic of $A(p, n, n-\alpha)$. Using~\eqref{eqn:masterEq}, we essentially show that for $k=n-\alpha$ we have
\[ \frac{n!}{(n-2\alpha)!}\frac{1}{\alpha!}\round{\frac{2^p-1}{2}}^{\alpha}\leq A(p, n, k)\leq \frac{n!}{(n-2\alpha)!}\frac{1}{\alpha!}\round{\frac{2^p-1}{2}}^{\alpha} + O_{p, \alpha}\round{n^{2\alpha-1}}\;.\]
We say a function $g(p, n, \alpha)=O_{p, \alpha}(f(n))$ if $g(p, n, \alpha\leq Cf(n)$ for some constant $C$ (possibly depending on $(p, \alpha)$). Because of this asymptotic, see that the first order term in $A(p, n, k)^2$ and $A(p, n, k-1)A(p, n, k+1)$ exactly agree (up to an explicit factor). More precisely, we define 
\[ M(p, n, \alpha)\coloneqq \frac{1}{\alpha!} \frac{n!}{(n-2\alpha)!}\round{\frac{2^p-1}{2}}^{\alpha}\;\;.\]
Note that $M(p, n, \alpha)= O_{p, \alpha}(n^{2\alpha})$ for fixed $p$ and $\alpha$. We essentially show that
\begin{align*}
    A(p, n, k)^2 &= M(p, n, \alpha)^2 + O_{p, \alpha}(n^{4\alpha-1})\\
    A(p, n, k-1)A(p, n, k+1) &=  \frac{\alpha}{\alpha+1}\frac{(n-2\alpha-1)(n-2\alpha)}{(n-2\alpha+2)(n-2\alpha+1)} M(p, n, \alpha)^2 + O_{p, \alpha}(n^{4\alpha-1})
\end{align*}

This clearly implies Theorem~\ref{thm:MainTheorem}. In this paper, we do not make any attempt to optimize the value of $N$. Also, the value of $N$, unfortunately, is not explicit yet, but we believe that our method can be refined to yield a concrete estimate of $N=N(\alpha, p)$. However, we should point out that the key issue in obtaining an estimate for $N(\alpha, p)$ comes from the fact that the hidden constants in the $O_{p, \alpha}(\cdot)$ depend on the function 
\[f(p, \alpha)\coloneqq \max\left\{\prod_{i=1}^{\alpha-1}\frac{B(p, \alpha_i)}{\alpha_i}: 1\leq \alpha_i\leq \alpha+1, \sum_{i=1}^{\alpha-1}\alpha_i=2\alpha-1\right\}\;.\]
While $\alpha$ is kept constant, this quantity is a constant and therefore we can argue that asymptotically $O_{p, \alpha}(n^{4\alpha-1})$ becomes negligible compared to the leading order term which is $O(n^{4\alpha})$. Of course, one can use a naive upper bound on $f(p, \alpha)$ to obtain some estimate for $N(\alpha, p)$; but we do not explore this direction here. Obtaining a good upper bound on $f(p, \alpha)$ is not only necessary for a reasonable estimate of $N(\alpha, p)$, but such investigations will be crucial to extend Theorem~\ref{thm:MainTheorem} to the regimes where $\alpha$ grows with $n$, say, $\alpha\approx \log(n)$ or $\alpha\approx \sqrt{n}$. We leave such investigations for future work.

\section{The case~\texorpdfstring{ $k=n-1,n-2$}{k=n-1, n-2}}
In this section, we prove Theorem~\ref{thm:MainTheorem} for the case $k=n-1$ and $k=n-2$. The Proof of the case $k=n-1$ is taken directly from~\cite{abdesselam2024bijection} and the proof for $k=n-2$ follows the same strategy and relies on exactly computing the $A(p, n, n-j)$ for $j=0, 1, 2, 3$. However, towards the end, we use the computations in these ``easy" cases to motivate a heuristic that inspires our proof of Theorem~\ref{thm:MainTheorem}. We skip the tedious computation for the case $k=n-1$ as it can be found in~\cite{abdesselam2024bijection}, but we give the full detail for the $k=n-2$ case for the convenience of the reader. Only Lemma~\ref{lem:LowerBound} and Lemma~\ref{lem:Bounds} proved in Section~\ref{subsec:Observations} are used in the proof of our main result. The reader can directly skip to~\ref{subsec:Observations} without affecting the understanding of the paper. However, we feel that the computations below are somewhat illuminating.

\subsection{The proof of the case~\texorpdfstring{$k=n-1, n-2$}{k=n-1, n-2}}
We begin by noting the following identities that can be easily verified from the definition of $B(p, m)$ given above
\begin{align*}
    B(p, 1)= 1\;,\qquad B(p, 2) =(2^p-1)\;,\qquad B(p, 3) = \frac{1}{2}(3^p-1)\;,\qquad B(p, 4)=\frac{1}{3}(2^p-1)(2^{p+1}-1)\;.
\end{align*}

Note recall from~\eqref{eqn:masterEq} that
\[ A(p, n, k) = \frac{n!}{u!} \times \sum_{n_1, \ldots, n_u\geq 1} \mathbbm{1}\{n_1+\ldots+n_u=n\}\times \prod_{i=1}^{u} \frac{B(p, n_i)}{n_i}\;.\]

In the following, we sometimes refer to a vector $(n_1, \ldots, n_k)$ as a configuration and we think of it as a configuration of balls into $k$ many distinguishable bins. Observe that any configuration $n_1+\ldots+n_u=n$ and $n_i\geq 1$ corresponds to putting $n$ indistinguishable balls into $u$ bins such that each bin is non-empty. To do so, we first begin by putting one ball in each bin and now the problem reduces to putting $n-u$ indistinguishable balls into $u$ bins. In the case when $u=n$, there is nothing to do, and all $n_i=1$. Therefore, $A(p, n, n)=1$. When $u=n-1$, after putting $1$ ball into each $u$ bin, we are left with $1$ ball that we can place in any of the $u$ bins. Therefore, there is exactly $1$ index $i\in [n-1]$ such that $n_i=1$ and the remaining $n_j$s are all $1$. In this case the product $\prod_{i=1}^{u}\frac{B(p, n_i)}{n_i}= \frac{B(p, 2)}{2}$. Since the special index $i$ such that $n_i=2$ can be chosen in ${n-1}$ ways. We obtain 
\begin{align}
    \label{eqn:Apnminus1}
     A(p, n, n-1) &= \frac{n!}{(n-1)!}{n-1\choose 1}\frac{(2^p-1)}{2}= {n\choose 2}(2^p-1)\;.
\end{align}

Now consider the case when $u=n-2$, in this case, we are left with putting $2$ balls into $n-2$ bins. There are two possibilities to this end. Firstly, we can find a bin (in $n-2$ ways) and put the remaining two balls into this bin. This gives the vectors $(n_1, \ldots, n_u)$ where all $n_i$ except one are equal to $1$ and the exceptional index is precisely $3$. This contributed $(n-2)\frac{B(p, 3)}{3!}=(n-2)\frac{3^p-1}{3!}$ in $A(p, n, n-2)$. The other possibility is to put the remaining two balls into two separate bins. This can be done in ${n-2\choose 2}$ ways and each such configuration yields exactly two $n_i$s that are equal to $2$ and the remaining $n_j$s are $1$. This kind of configuration, therefore, contributed ${n-2\choose 2}\round{\frac{B(p, 2)}{2}}^2$. Thus, we obtain   
\begin{align}
\label{eqn:Apnminus2}
    A(p, n, n-2) &= \frac{n!}{(n-2)!}\round{{n-2\choose 1}\frac{(3^p-1)}{3!}+ {n-2\choose 2} \round{\frac{(2^p-1)}{2}}^2}\nonumber\\
    &= {n\choose 3}(3^p-1)+3{n\choose 4}(2^p-1)^2\;.
\end{align}
Using~\eqref{eqn:Apnminus1} and~\eqref{eqn:Apnminus2} and a straightforward computation that is performed in~\cite[Proposition 3.1]{abdesselam2024bijection}, one can verify the Conjecture~\ref{conj:Conjecture} for $k=n-1$. We leave the details to the interested reader.

We will give more details about the computations for the case $k=n-2$ for the convenience of the reader as it did not appear previously. To begin, we must compute $A(p, n, n-3)$. Once again following our previous analogy, after putting $1$ balls into each of $n-3$ bins we are left with the $3$ balls to be put into $n-3$ bins. There are three classes of possible configurations for this:
\begin{enumerate}
    \item One can distribute these three balls into three distinct bins. In this case, we get $(n_1, \ldots, n_{n-3})$ with exactly three $2$ and remaining $1$s. These three indices can be chosen in ${n-3\choose 3}$ ways. For each configuration in this class we have $\prod_{i=1}^{n-3}\frac{B(p, n_i)}{n_i}=\round{\frac{B(p, 2)}{2}}^3=\round{\frac{2^p-1}{2}}^3$. 

    \item One can put $1$ ball into one bin and $2$ balls into another. This gives $n_1, \ldots, n_k$ with exactly one position being $2$ and $3$ each, and the remaining positions are $1$. This can be done in $2\times {n-3\choose 2}$ ways. For any such configuration, we have $\prod_{i=1}^{n-3}\frac{B(p, n_i)}{n_i} =\frac{B(p, 2)}{2}\frac{B(p, 3)}{3}=\frac{(2^p-1)(3^p-1)}{2\times 2\times 3}$.

    \item Finally, one can put all the remaining $3$ balls into a single bin. This yields configurations where all but one position is $1$ and there is exactly one position with $n_i=4$. There are $(n-3)$ such configurations and for each such configuration we have $\prod_{i=1}^{n-3}\frac{B(p, n_i)}{n_i} = \frac{B(p, 4)}{4}= \frac{(2^p-1)(2^{p+1}-1)}{3\times 4}$.
\end{enumerate}


Summing over all the above configurations, we thus obtain 
\begin{align}
\label{eqn:Apnminus3}
    A(p, n, n-3) &= \frac{n!}{(n-3)!}\round{{n-3\choose 3}\round{\frac{2^p-1}{2}}^3+ 2{n-3\choose 2}\frac{(2^p-1)(3^p-1)}{2\times 2\times 3}+{n-3\choose 1}\frac{(2^p-1)(2^{p+1}-1)}{3\times 4}}\nonumber\\
    &= 2{n\choose 4}(2^{p}-1)+4{n\choose 4}(2^p-1)^2+ 10{n\choose 5}(2^p-1)(3^p-1)+15{n\choose 6}(2^p-1)^3\;.
\end{align}

Using~\eqref{eqn:Apnminus1},~\eqref{eqn:Apnminus2} and~\eqref{eqn:Apnminus3}, we compute
\begin{align*}
    A(p, n , n-2)^2 &= {n\choose 3}^2(3^p-1)^2+6{n\choose 3}{n\choose 4}(2^p-1)^2(3^p-1)+9{n\choose 4}^2(2^p-1)^4\;,\\
    A(p, n, n-1)A(p, n, n-3)&= 2{n\choose 2}{n\choose 4}(2^{p}-1)^2+4{n\choose 2}{n\choose 4}(2^p-1)^3\\
    &\qquad + 10{n\choose 2}{n\choose 5}(2^p-1)^2(3^p-1)+15{n\choose 2}{n\choose 6}(2^p-1)^4\;.
\end{align*}

Now observe that
\begin{align}
\label{eqn:Delta}
    \Delta(p, n, n-2) &= A(p,n, n-2)^2-A(p, n, n-1)A(p, n, n-3)\\
    &= \round{9{n\choose 4}^2 - 15{n\choose 2}{n\choose 6}}(2^p-1)^4\nonumber\\
    &+\round{6{n\choose 3}{n\choose 4}- 10{n\choose 2}{n\choose 5}}(3^p-1)(2^p-1)^2\nonumber\\
    &+\round{{n\choose 3}^2(3^p-1)^2-2{n\choose 2}{n\choose 4}(2^p-1)^2-4{n\choose 2}{n\choose 4}(2^p-1)^3}\nonumber\;.
\end{align}

It is easily verified that
\begin{align*}
   A(n)&\coloneqq \round{9{n\choose 4}^2 - 15{n\choose 2}{n\choose 6}}\\
   &= \frac{1}{2}{n\choose 2}{n\choose 4}\round{\frac{3}{2}(n-2)(n-3)-(n-4)(n-5)}\\
   &=\frac{1}{2}{n\choose 2}{n\choose 4}\left(\frac{n^2}{2}+\frac{3n}{2}-11\right) = \frac{1}{4}(n^2+3n-22){n\choose 2}{n\choose 4}\;.
   \end{align*}
   Note that $n^2+3n-22$ is increasing in $n$, for all positive $n$, and for $n=4$, we get $n^2+3n-22=6$. In particular, $A(n)\geq \frac{6}{4}{n\choose 2}{n\choose 4}$. 
   Now, note that
   \begin{align*}
B(n) &\coloneqq \round{6{n\choose 3}{n\choose 4}- 10{n\choose 2}{n\choose 5}}\\
&=\frac{n^2(n-1)^2(n-2)(n-3)}{24}\round{(n-2)-(n-4)}\\
&=\frac{n^2(n-1)^2(n-2)(n-3)}{12}=4{n\choose 2}{n\choose 4}\;.
\end{align*}
Combining these with equation~\eqref{eqn:Delta}, we obtain 
\begin{align*}
 \Delta(p, n, n-2) &= A(n)(2^p-1)^4+B(n)(3^p-1)(2^p-1)^2\\
 &+\left({n\choose 3}^2(3^p-1)^2 - \frac{1}{2}B(n)(2^p-1)^2-B(n)(2^p-1)^3\right)\\
 &\geq A(n)(2^p-1)^4 +{n\choose 3}^2(3^p-1)^2 +B(n)(2^p-1)^2\round{(3^p-1)-\frac{1}{2}-(2^p-1)}\\
 &=A(n)(2^p-1)^4 +{n\choose 3}^2(3^p-1)^2 + B(n)(2^p-1)^2\round{3^p-2^p-\frac{1}{2}}\;.
\end{align*}

For $p\geq 2$, we have that $3^p-2^p-\frac{1}{2}\geq 0$, it follows that $\Delta(p, n, n-2)\geq 0$. Of course, we have a much stronger bound on $\Delta(p, n, n-2)$. In fact $\Delta(p, n, n-2)\sim n^6(2^p-1)^4$.

\subsection{Some observations and useful preliminaries}
\label{subsec:Observations}
Observe the expression for $A(p, n, n-\alpha)$ for $\alpha=1, 2, 3$ in~\eqref{eqn:Apnminus1},~\eqref{eqn:Apnminus2} and~\eqref{eqn:Apnminus3}. 
Note that the leading order term (in $n$) in $A(p, n, n-\alpha)$ is \[\frac{n!}{(n-\alpha)!}{n-\alpha\choose \alpha}\round{\frac{(2^p-1)}{2}}^\alpha\;.\] 
Where does this term really come from? This leading order term from the case when $n_i\in \{1, 2\}$ for all $i\in [k]$. In other words, there are exactly $\alpha$ many positions with $2$ and the remaining $n-2\alpha$ positions are $1$. Let us denote the collection of all such vectors/configurations by $\Ccal(n, \alpha)$. Since $\abs{\Ccal(n, \alpha)}={n-\alpha\choose \alpha}$ and each configuration in $\Ccal(n, \alpha)$ contributes precisely $\round{\frac{2^p-1}{2}}^\alpha$. This explains the leading order term. 

Let us denote all the remaining configurations $(n_1, \ldots, n_{n-\alpha})$ by $\Dcal(n, \alpha)$. Obtaining a given configuration $(n_1, \ldots, n_{n-\alpha})\in \Dcal$ is equivalent to first putting $1$ ball into each of $n-\alpha$ bins. Then, the remaining $\alpha$ balls have to be put into at most $(\alpha-1)$ positions.  This gives $\abs{\Dcal(n, \alpha)}\leq {n-\alpha\choose \alpha-1}{2\alpha-2\choose \alpha}$. Note that, in general, this inequality is not sharp and possibly can be improved. 

Our goal now is to control  $\prod_{i=1}^{n-\alpha}\frac{B(p, n_i)}{n_i}$
over all configurations $(n_1, \ldots, n_{n-\alpha})$ in $\Dcal(n, \alpha)$. To this end, fix a configuration $(n_1, \ldots, n_{n-\alpha})$ in $\Dcal(n, \alpha)$ and suppose $\alpha_1, \ldots, \alpha_t$ are the values that are distinct from $1$ where $1\leq t\leq \alpha-1$. For such configuration we have $\alpha_1+\ldots+\alpha_t=(\alpha+t)$ and $2\leq \alpha_i\leq \alpha+1$. For such a configuration, we have 
\begin{equation}
\label{eqn:InnerProd}
    \prod_{i=1}^{n-\alpha}\frac{B(p, n_i)}{n_i}= \prod_{i=1}^{t}\frac{B(p, \alpha_i)}{\alpha_i}\;.
\end{equation}

To bound this, we define
\begin{equation}
    \label{eqn:UpperBoundMultipicativeFun}
    f(p, \alpha)\coloneqq \max\left\{\prod_{i=1}^{\alpha-1}\frac{B(p, \alpha_i)}{\alpha_i}: 1\leq \alpha_i\leq \alpha+1, \sum_{i=1}^{\alpha-1}\alpha_i=2\alpha-1\right\}\;.
\end{equation}
With this notation, we obtain that for any configuration in $\Dcal(n, \alpha)$, the product in~\eqref{eqn:InnerProd} is at most $f(p, \alpha)$. Note that $f(p, \alpha)$ is independent of $n$ for every fixed $p$ and $\alpha$. This is important for our argument. 

Combining the above computations, we obtain
\begin{align*}
    A(p, n, n-\alpha) &\geq \frac{n!}{(n-\alpha)!}{n-\alpha\choose \alpha}\round{\frac{(2^p-1)}{2}}^\alpha\;,\\
    A(p, n, n-\alpha) &\leq \frac{n!}{(n-\alpha)!}{n-\alpha\choose \alpha}\round{\frac{(2^p-1)}{2}}^\alpha+ \frac{n!}{(n-\alpha)!}{n-\alpha\choose \alpha-1}{2\alpha-2\choose \alpha} f(p, \alpha)\;.
\end{align*}

We record the above discussion as the following lemma that will be useful later. Since we have already proved the conjecture for the case $k=n-1$ and $k=n-2$, for the following lemma, we will work with $k=n-\alpha$ for some $\alpha\geq 3$.

\begin{lemma}
\label{lem:LowerBound}
Fix $\alpha\geq 3$ and $p\geq 2, n\geq 4$. Set $k=n-\alpha$. Then, we have
    \[A(p, n, k)\geq \frac{n!}{\alpha!(n-2\alpha)!}\round{\frac{2^p-1}{2}}^\alpha\;.\]
\end{lemma}

\begin{lemma}
\label{lem:Bounds}
    Fix $\alpha\geq 3$ and $p\geq 2, n\geq 2\alpha$. Set $k=n-\alpha$. Then, we have
    \begin{align*}
        A(p, n, k-1)&\leq \frac{n!}{(\alpha+1)!(n-2\alpha-2)!}\round{\frac{2^p-1}{2}}^{\alpha+1}+ \frac{n!}{(n-2\alpha-1)!\alpha!}{2\alpha\choose \alpha+1} f(p, \alpha+1)\;,\\
        A(p, n, k+1)&\leq \frac{n!}{(\alpha-1)!(n-2\alpha+2)!}\round{\frac{2^p-1}{2}}^{\alpha-1}+ \frac{n!}{(n-2\alpha+3)!(\alpha-2)!}{2\alpha-4\choose \alpha-1}f(p, \alpha-1)\;.
    \end{align*}
\end{lemma}
The proofs follow from the previous discussion so we skip the details. 

\section{Proof of Theorem~\ref{thm:MainTheorem}}
In this section, we prove Theorem~\ref{thm:MainTheorem}. Our proof uses Lemma~\ref{lem:LowerBound} and Lemma~\ref{lem:Bounds} from Section~\ref{subsec:Observations}.

We begin by recalling the setting and establishing some notations for simplicity. 
Recall that $\alpha\geq 3, p\geq 2$ are fixed and we set $k=n-\alpha$. Now set 
\[ M(p, n, \alpha) \coloneqq \frac{n!}{\alpha!(n-2\alpha)!}\round{\frac{2^p-1}{2}}^\alpha =\frac{(n)_{2\alpha}}{\alpha!}\round{\frac{2^p-1}{2}}^\alpha \;,\]
where we used the Pochhammer symbol $(n)_{x}$ to denote $(n)_{x}=n(n-1)\cdots (n-x+1)$. We note for fixed $x$, asymptotically, $(n)_x\sim n^{x}$, that is, $\frac{(n)_{x}}{n^x}\to 1$ as $n\to \infty$. 

With this notation, we re-write the bounds from the Lemma~\ref{lem:LowerBound} and Lemma~\ref{lem:Bounds} as follows:
 \begin{align*}
       A(p, n, k)&\geq M(p, n, \alpha)\;,\\
        A(p, n, k-1)&\leq \frac{(n-2\alpha)_2}{(\alpha+1)}\round{\frac{2^p-1}{2}} M(p, n, \alpha)+ \frac{(n)_{2\alpha+1}}{\alpha!}{2\alpha\choose \alpha+1} f(p, \alpha+1)\;,\\
        A(p, n, k+1)&\leq \frac{\alpha}{(n-2\alpha+2)_2}\round{\frac{2}{2^p-1}} M(p, n, \alpha)+\frac{(n)_{2\alpha-3}}{(\alpha-2)!}{2\alpha-4\choose \alpha-1}f(p, \alpha-1)\;.
    \end{align*}

We set up some more notations for convenience.
\begin{align*}
    X &\coloneqq  \frac{\alpha}{\alpha+1}\frac{ (n-2\alpha)_2}{(n-2\alpha+2)_2}\;, \\
    Y_{1} &\coloneqq (n)_{2\alpha-3}(n-2\alpha)_2\times \left[\frac{\alpha(\alpha-1)}{(\alpha+1)!}{2\alpha-4\choose \alpha-1}\round{\frac{2^p-1}{2}}f(p, \alpha-1)\right]\;,
     \\
    Y_2 &\coloneqq \frac{(n)_{2\alpha+1}}{(n-2\alpha+2)_2}\times \left[\frac{1}{(\alpha-1)!}{2\alpha\choose \alpha+1}\round{\frac{2}{2^p-1}}f(p, \alpha+1)\right]\;,\\
    Z &\coloneqq (n)_{2\alpha+1}(n)_{2\alpha-3}\times\left[\frac{1}{\alpha!(\alpha-2)!}{2\alpha\choose \alpha+1}{2\alpha-4\choose \alpha-1}f(p, \alpha-1)f(p, \alpha+1)\right]\;.
\end{align*}

With this setup, we have 
\begin{align*}
    \frac{A(p, n, k-1)A(p, n, k+1)}{A(p, n, k)^2}&\leq \frac{XM(p, n, \alpha)^2+ (Y_1+Y_2)M(p, n, \alpha)+ Z}{M(p, n, \alpha)^2}\\
    &= X+ \frac{Y_1+Y_2}{M(p, n, \alpha)}+ \frac{Z}{M(p, n, \alpha)^2}\;.
\end{align*}

Now notice that 
\begin{align*}
    X =\frac{\alpha}{\alpha+1}\round{1-\frac{2}{n-2\alpha+1}}\round{1-\frac{2}{n-2\alpha+2}}<\frac{\alpha}{1+\alpha}\;,
\end{align*}
and as $n\to \infty$ we have that $X\to \frac{\alpha}{1+\alpha}$. 

Now observe that for fixed $p$ and $\alpha$, the variable $Y_1$ and $Y_2$ are of order $n^{2\alpha-1}$ asymptotically in $n$, while $M(p, n, \alpha)$ is of order $n^{2\alpha}$. It follows that 
\[ \frac{Y_1+Y_2}{M(p, n, \alpha)}\to 0,\]
as $n\to \infty$. 

Similarly, we also note that $Z$ is of order $n^{4\alpha-2}$ while $M^2$ is of order $n^{4\alpha}$. Since $\alpha$ and $p$ are kept fixed, we observe that $\frac{Z}{M(p, n, \alpha)}\to 0$ as $n\to \infty$.


In particular, we conclude that for every fixed $p\geq 2$ and $\alpha\geq 3$, we have 

\[\lim_{n\to \infty} \frac{A(p, n, k-1)A(p, n, k+1)}{A(p, n, k)^2}\leq \frac{\alpha}{1+\alpha}<1.\]

Since the limit above exists and is $<1$, it follows that there exists a finite $N=N(\alpha, p)$ such that $\frac{A(p, n, k)^2}{A(p, n, k-1)A(p, n, k+1)}\geq 1$ for $n\geq N$.

\newpage


\bibliographystyle{alpha} 
\bibliography{references}

\end{document}